\documentclass[a4paper,twoside,11pt,openright]{article}

\usepackage{bea_article}

\begin{document}

\title{Partition function of periodic isoradial dimer models}
\author{B\'eatrice de Tili\`ere
\thanks{Supported by Swiss National Fund under grant 47102009}\vspace{0.4cm}\\
       {\small Institut f\"ur Mathematik,
       Universit\"at Z\"urich,
       Winterthurerstrasse 190,
       CH-8057 Z\"urich.}\\
       {\small\texttt{beatrice.detiliere@math.unizh.ch}}
}
\date{}
\maketitle

\vspace{-1cm}
\begin{abstract}
Isoradial dimer models were introduced in \cite{Kenyon3} - they
consist of dimer models whose underlying graph satisfies a simple
geometric condition, and whose weight function is chosen
accordingly. In this paper, we prove a conjecture of \cite{Kenyon3},
namely that for periodic isoradial dimer models, the growth
rate of the toroidal partition function has a simple explicit formula
involving the local geometry of the graph only. This is a surprising
feature of
periodic isoradial dimer models, which does not hold in the general
periodic dimer case \cite{KOS}.
\end{abstract}

\section{Introduction}\label{sec0}

\noindent In this paper, we solve a conjecture of \cite{Kenyon3},
namely that the growth rate of the partition function of
periodic isoradial dimer models can be expressed explicitly, using only
the local geometry of the underlying graph. This is a surprising
feature of isoradial dimer models, which does not hold in the case of
general periodic dimer models. Indeed, in the general periodic case,
Kenyon, Okounkov and Sheffield
\cite{KOS} obtain an expression for the growth rate of the partition
function involving elliptic integrals, thus making it hard to do
explicit computations. An example of application of our
result yields a two line proof of Kasteleyn's celebrated result for
the growth rate of the dimer partition function of the
quadratic lattice $\ZZ^2$ \cite{Kasteleyn1, Kasteleyn2}, see Section
\ref{subsec05}. In order to state our result, let us precisely
describe the setting.

\subsection{Dimer model}\label{subsec01}

\begin{figure}[ht]
\begin{center}
\includegraphics[width=9cm]{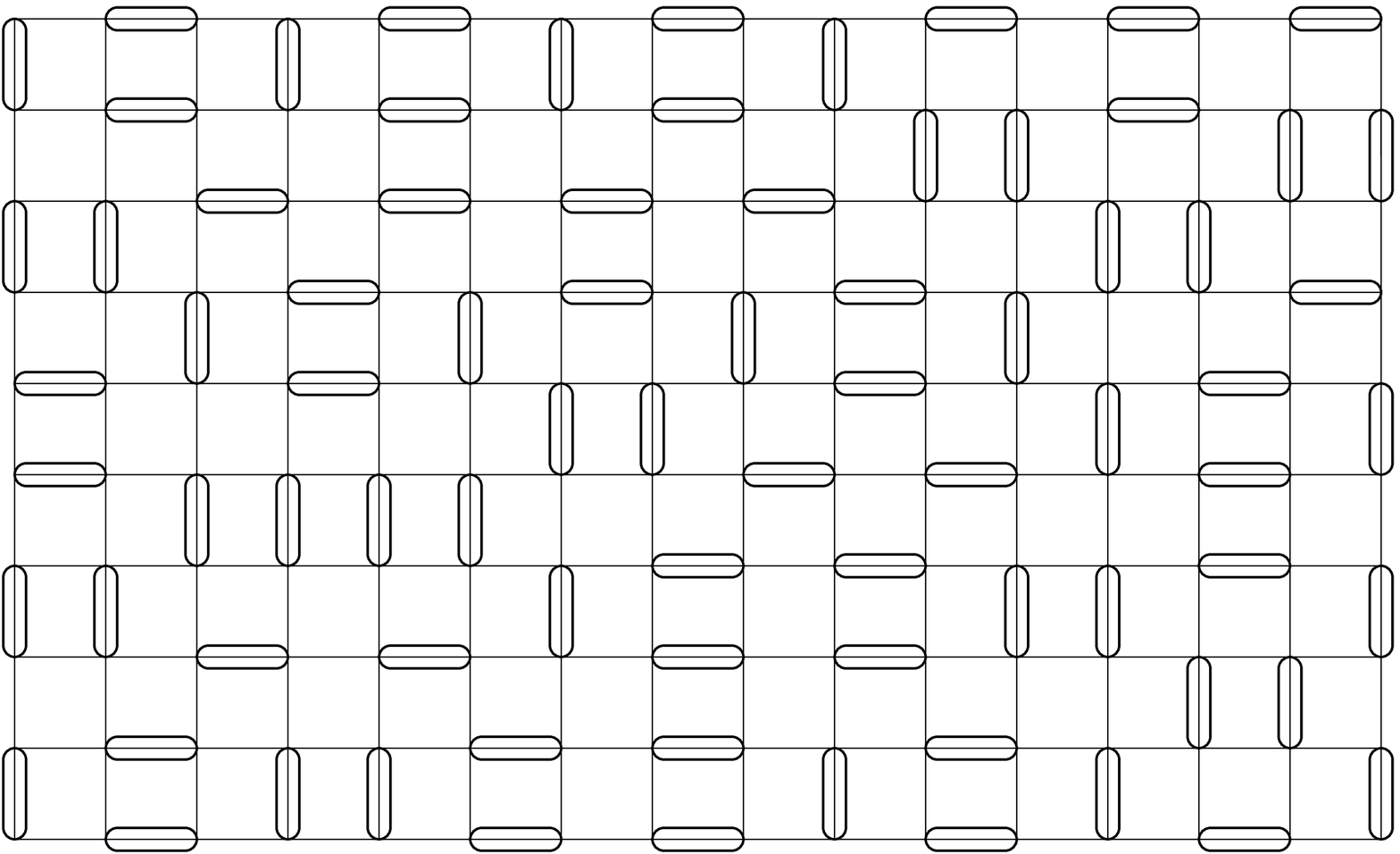}
\end{center}
\caption{A dimer configuration of a finite subgraph of $\ZZ^2$.} \label{fig2}
\end{figure}

\noindent The {\bf dimer model} belongs to the field of statistical
mechanics, and represents the adsorption of diatomic molecules on the
surface of a crystal, it is defined in the following way. 
The surface of the crystal is modeled by a graph $G$. We assume that
$G$ is simple, infinite, simply connected (i.e. it is the one skeleton of a simply
connected union of faces) and that its vertices are of degree
$\geq 3$. A {\bf dimer configuration} of $G$ is a {\bf perfect
matching} of $G$, that is a subset of edges $M$, such that every
vertex of $G$ touches a unique edge of $M$; refer to Figure \ref{fig2}
for an example in the case where $G$ is a finite subgraph of $\ZZ^2$. Let
us denote by $\M(G)$ the set of perfect matchings of $G$.\jump

\noindent Consider a simply connected finite subgraph $G_0$ of $G$. Then, dimer configurations of $G_0$ are
chosen with respect to the {\bf Boltzmann measure}, defined in the
following way.
Assume that a positive weight function $\nu$ is assigned to
edges of $G$, that is, each edge $e$ of $G$ has a weight
$\nu(e)$. Then, each dimer configuration
$M$ of $G_0$ has an energy $\E(M)=-\sum_{e\in M} \log \nu(e)$. The
probability of occurrence of the dimer configuration $M$ chosen with
respect to the Boltzmann measure $\mu_0$, is given by:
\begin{equation*}
\mu_0(M)=\frac{e^{-\E(M)}}{Z(G_0,\nu)}=\frac{\prod_{e\in M}\nu(e)}{Z(G_0,\nu)},
\end{equation*}
where $Z(G_0,\nu)=\sum_{M\in\M(G_0)}\prod_{e\in M}\nu(e)$ is the
normalizing constant, known as the {\bf partition function}.\jump
The first step in the study of the dimer model is the computation of the
partition function. Indeed, the latter yields precious information
about the global behavior of the system \cite{Baxter,CKP,KOS}, whose understanding is the
goal of statistical mechanics.
Kasteleyn \cite{Kasteleyn1, Kasteleyn2}, and independently Temperley \& Fisher \cite{TF},
laid the ground stone for
the study of the partition function. Consider $G_0$ as above. Then
they give an explicit formula for $Z(G_0,\nu)$ as the square root of
the determinant of a matrix $\Ks_0$, also known as a {\bf Kasteleyn
matrix}. Loosely stated, $\Ks_0$ is the weighted adjacency matrix of
$G_0$ (weighted by the function $\nu$), where minus signs are added to
coefficients in a suitable way, see Section \ref{subsec12} for
details.\jump
Since edges of dimer configurations represent diatomic molecules, the
next step is to understand the partition function of infinite
graphs. The relevant question then, is the computation of its growth
rate. More precisely, if $\{G_n\}$
is an exhaustion of $G$ by finite graphs, the goal is to compute: 
\begin{equation}\label{equ1}
C=\lim_{n\rightarrow\infty}\frac{1}{n^2}\log Z(G_n,\nu).
\end{equation}
Note that this limit is not
universal, in the sense that it strongly depends on the choice of
exhaustion $\{G_n\}$. Moreover, in the case of exhaustions by planar
graphs, this limit is hard to compute because of
the lack of symmetry which reflects in the Kasteleyn matrix. This is
the reason why computations are done on toroidal graphs, as explained
in the next section.

\subsection{Toroidal dimer models}\label{subsec02}

\noindent Consider an infinite graph $G$ as above, and suppose
that $G$ is bipartite, i.e that it admits a bipartite coloring of its
vertices. Assume moreover that there exists a bi-dimensional lattice
$\Lambda$, such that $G$ is doubly $\Lambda$-periodic, that is the graph 
$G$ as well as its vertex coloring are bipartite. Define $G_n$ to be
the toroidal
graph $G/\Lambda n$, then $\{G_n\}$ is a natural exhaustion of the
graph $G$.\jump
In the case where $G$ is the square lattice $\ZZ^2$, Kasteleyn
\cite{Kasteleyn1,Kasteleyn2} gives an explicit expression for
$Z(G_n,\nu)$ as a linear combination of the determinants of four
Kasteleyn matrices. Tesler \cite{Tesler} generalizes this result to
any graph $G$ satisfying the above, see Section \ref{subsec12} for
details (he actually generalizes it to graphs embedded on genus $g$ 
surfaces). In \cite{KOS}, Kenyon, Okounkov and Sheffield give an
explicit expression for the limit $C$ of equation (\ref{equ1}) in the
case of the natural toroidal exhaustion $\{G_n\}$ of $G$. For
general periodic, bipartite graphs $G$, this limit depends on the
combinatorics of the model and involves 
elliptic integrals.\jump
In the next section, we define isoradial dimer models - the sub-family
of dimer models for which the expression for the growth rate $C$ of (\ref{equ1})
becomes surprisingly simple.

\subsection{Isoradial dimer models}\label{subsec03}

\noindent Isoradial dimer models were introduced in
\cite{Kenyon3}. They are dimer models on graphs $G$ satisfying a geometric
condition called {\bf isoradiality}. This notion first appeared in
\cite{Duffin}, see also \cite{Mercat0}, and is defined as follows. All faces of $G$ are
inscribable in a circle, and all circumcircles have the same radius,
moreover, all circumcenters of the faces are contained in the closure
of the faces. The common radius is taken to be $1$. An isoradial 
embedding of the dual graph $G^*$ is obtained by
sending dual vertices to the center of the corresponding faces.\jump
Recall that the energy of a dimer configuration depends on the weights
assigned to edges of $G$. In the case of isoradial graphs,
one considers a specific weight function $\nu$, called the {\bf critical
weight function}, defined as follows \cite{Kenyon3}. To each edge $e$
of $G$ corresponds a unit side-length rhombus $R(e)$ whose vertices are the
vertices of $e$ and of its dual edge ($R(e)$ may be degenerate). Let
$\widetilde{R}=\cup_{e\in G}R(e)$. Then, define $\nu(e)=2\sin\theta$,
where $2\theta$ is the angle of the rhombus $R(e)$ at the vertex it
has in common with $e$; $\theta$ is called the {\bf rhombus angle} of
the edge $e$. Note that $\nu(e)$ is the length of $e^*$ the dual edge
of $e$.\jump
When the graph considered for a dimer model is isoradial, and when the
energy of configurations is determined by the critical weight
function, we speak of an {\bf isoradial dimer model}.

\subsection{Result}\label{subsec04}

Consider an isoradial dimer model. Assume moreover that $G$ is
bipartite and doubly $\Lambda$-periodic. As above consider the natural
exhaustion $\{G_n\}$ of $G$ by toroidal graphs. Then, Theorem
\ref{thm1} below gives an explicit expression for the growth rate of the
partition function of the exhaustion $\{G_n\}$ as a function of the
local geometry of the graph.
\begin{thm}\label{thm1}
\begin{equation}\label{equ2}
\lim_{n\rightarrow\infty}\frac{1}{n^2}\log Z(G_n,\nu)=
\sum_{i=1}^m\left(\frac{\theta_i}{\pi}\log
  2\sin\theta_i+\frac{1}{\pi}L(\theta_i)\right),
\end{equation}
where, $\theta_1,\ldots,\theta_m$ are the rhombus angles of the edges
$e_1,\ldots,e_m$ of $G_1$, and $L$ is Lobachevsky's function,
$L(x)=-\int_{0}^x\log 2\sin t dt$.
\end{thm}
\begin{itemize}
\item The right hand side of (\ref{equ2}) is (up to a factor $2$), the expression Kenyon obtained
for what he calls the {\em log of the normalized determinant} of the
Dirac operator \cite{Kenyon3}. He conjectured it to be the right limit for
the growth rate of the partition function of an appropriate
exhaustion of $G$. Theorem \ref{thm1} states that this is the right limit
in the case of the natural exhaustion $\{G_n\}$ of $G$ by toroidal graphs. 
\item The surprising feature of Theorem \ref{thm1} is that the growth
rate of the partition function is expressed using geometric
information of the fundamental domain $G_1$ only, no combinatorics
is involved. This is contrast with the result of \cite{KOS} for
general bipartite periodic graphs. Note that, although the
result of \cite{KOS} is true in a more general setting, there seems
to be no simple way of re-deriving Theorem \ref{thm1} from it. Let
us also mention that we use the result of \cite{KOS} in order to have a priori existence of the limit
$C$ of equation (\ref{equ1}).
\end{itemize}

\subsection{Example: Kasteleyn's computation}\label{subsec05}

\begin{figure}[ht]
\begin{center}
\includegraphics[width=10cm]{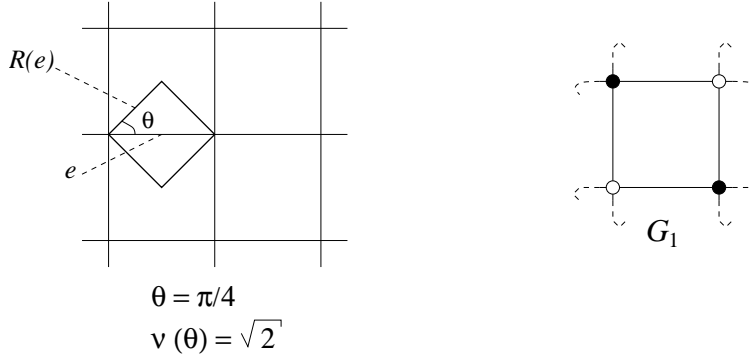}
\end{center}
\caption{Critical weight function for $\ZZ^2$ (left). Fundamental
domain $G_1$ (right).} \label{fig1}
\end{figure}

\noindent In \cite{Kasteleyn1,Kasteleyn2}, Kasteleyn skillfully
computes the toroidal partition function for $G=\ZZ^2$. As an example
of application of Theorem \ref{thm1}, let us rederive his result in a
two line calculation. Note
that the following computation was already mentioned in \cite{Kenyon3},
as a support for the conjecture proved in Theorem \ref{thm1}. 
When $G=\ZZ^2$, the critical weight function is $\nu\equiv\sqrt{2}$, see
Figure \ref{fig1} (left). Moreover, consider the fundamental domain
$G_1=G/\Lambda$ as in Figure \ref{fig1} (right). 

\noindent Then, a direct
computation using Theorem \ref{thm1} yields:
\begin{eqnarray*}
\lim_{n\rightarrow\infty}\frac{1}{2n^2}\log
Z(G_n,\sqrt{2})&=&4\left(\frac{1}{4}\log\sqrt{2}+\frac{1}{\pi}L\left(\frac{\pi}{4}\right)\right),\\
&=&\frac{1}{2}\log 2 +\frac{2\cal K}{\pi},
\end{eqnarray*}
where $\cal K$ is Catalan's constant. Note that the $\frac{1}{2}\log 2$ factor differs from
Kasteleyn's computation, because we have weights $\sqrt{2}$ on the
edges instead of weights $1$.\jump
{\em Acknowledgments}: we would like to thank Richard Kenyon
for proposing the question solved in this paper, and for fruitful discussions.

\section{Proof of Theorem \ref{thm1}}\label{sec1}

\noindent In the whole of this section, we let $G$ be an infinite, bipartite,
$\Lambda$-periodic, isoradial graph. We assume that edges of $G$ are
assigned the critical weight function. Moreover, $W$ denotes the set of white
vertices of $G$, $B$ the set of black ones, and $G_n$ the
toroidal graph $G/\Lambda n$.\jump
Before giving the proof of Theorem \ref{thm1}, let us
precisely describe the conjecture of \cite{Kenyon3}. In \cite{Kenyon3},
Kenyon introduces the Dirac operator $K$ indexed by vertices of $G$
(we refer the reader to \cite{Kenyon3} for the definition of $K$, and
to Section \ref{subsec11} for the definition of the {\em real Dirac
operator} $\Ks$ which is closely related).
He also defines its inverse $K^{-1}$, for which he proves existence
and uniqueness. One of his beautiful result is an explicit expression
for the {\em log of
the normalized determinant} ``log det$_1$ $K$'' of the Dirac operator
$K$, defined by the following PDE:
\begin{equation}\label{equ3}
\frac{\partial \log \det_1 K}{\partial K(w,b)}=\frac{1}{|V(G_1)|}K^{-1}(b,w),
\end{equation}
where $w\in W$, $b\in B$ are adjacent vertices, $|V(G_1)|$ is the
number of vertices of $G_1$. Adding some initial condition allows him
to solve the PDE, and obtain the explicit expression for $\log \det_1
K$, which is the right hand side of equation (\ref{equ2}). How is
this related to the dimer model? The Dirac operator $K$ is closely related
to an infinite Kasteleyn matrix - to have an actual infinite Kasteleyn
matrix, one needs to work with the {\em real Dirac operator} $\Ks$
(see Section \ref{subsec11} below), obtained from the Dirac operator by a
gauge transformation. If $G_0$ is a finite subgraph of $G$, and
$\Ks_0$ the restriction of $\Ks$ to the vertices of $G_0$, then
$\log\det_1 \Ks_0$ can be defined naturally as
$\frac{1}{|V(G_0)|}\log\det \Ks_0$, which is the normalized log of the
dimer partition function of $G_0$. Hence, the conjecture is to
interpret (\ref{equ3}) as the limit of the dimer partition function on some
appropriate exhaustion of $G$. The problem lies in proving the
existence of such limits, leading us to work with toroidal
exhaustions. Note that part of
the proof is close to \cite{Kenyon3}, since it consists in solving the
PDE (\ref{equ3}) for operators
restricted to subgraphs of the exhaustion.\jump
The structure of Section \ref{sec1} is as follows. 
Section \ref{subsec11} consists in the definition of the real Dirac
operator $\Ks$ and its inverse $\Ks^{-1}$. The operator $\Ks$ was introduced in \cite{Bea}, see also
\cite{KO}. Note that features of $\Ks$ are closely related to those of the Dirac
operator $K$ of \cite{Kenyon3}. In Section \ref{subsec12}, we state the
explicit expression of \cite{Kasteleyn2, Tesler} for the toroidal
partition function. Section \ref{subsec13} uses Sections
\ref{subsec11} and \ref{subsec12}, and consists in the
proof of Theorem \ref{thm1}.

\subsection{Real Dirac operator}\label{subsec11}

\noindent The {\em real Dirac operator} $\Ks$ is obtained from the {\em
Dirac operator} $K$ of \cite{Kenyon3} by a gauge
transformation. Both operators are represented by weighted infinite
adjacency matrices indexed by vertices of $G$. For $K$ the edges are
unoriented and weighted by the critical weight function times a
complex number of modulus $1$. For $\Ks$, edges are oriented with a
{\em clockwise odd} orientation (see below), and are weighted by the critical
weight function. Both weight functions yield the same
Boltzmann measure on finite simply connected sub-graphs of $G$, but,
and this is the reason why the real Dirac operator $\Ks$ is introduced,
these weights do not yield the same probability distribution on
toroidal subgraphs. Note that the real Dirac operator was already used
in \cite{Bea}, see also \cite{KO}. 

\subsubsection{Definition}\label{subsec111}

\noindent Following Kasteleyn \cite{Kasteleyn2}, let us define
clockwise odd orientations on edges of $G$. An
elementary cycle $\C$ of $G$ is said to be {\bf clockwise odd}
if, when traveling cw (clockwise) around the edges of $\C$, the number of
co-oriented edges is odd. Note that since $G$ is bipartite, the
number of contra-oriented edges is also odd. Then, an orientation of
the edges of $G$ is defined to be {\bf clockwise odd} if all 
elementary cycles are clockwise odd. In \cite{Kasteleyn2}, Kasteleyn
shows that, for planar simply connected graphs, such an
orientation always exists.\jump
Consider a clockwise odd orientation of the edges of $G$. Define
$\Ks$ to be the infinite adjacency matrix of the graph $G$,
weighted by the critical weight function $\nu$. That is, if $v_1$
and $v_2$ are not adjacent, $\Ks(v_1,v_2)=0$. If $w\in W$ and
$b\in B$ are adjacent vertices, then
$\Ks(w,b)=-\Ks(b,w)=(-1)^{\II_{(w,b)}}\nu(wb)$, where
$\II_{(w,b)}=0$ if the edge $wb$ is oriented from $w$ to $b$, and
$1$ if it is oriented from $b$ to $w$. The infinite matrix $\Ks$
defines the {\bf real Dirac operator} $\Ks$: $\CC^{V(G)}
\rightarrow \CC^{V(G)}$, by
\begin{equation*}\label{618}
(\Ks f)(v)=\sum_{u \in G} \Ks(v,u)f(u).
\end{equation*}
The matrix $\Ks$ is also called a {\bf Kasteleyn matrix} for the
underlying dimer model.

\subsubsection{Inverse real Dirac operator}\label{subsubsec112}

\noindent Similarly to the definition of the inverse Dirac operator $K^{-1}$
 \cite{Kenyon3}, the {\bf inverse real Dirac
 operator} $\Ks^{-1}$ is
defined to be the unique operator satisfying:
\begin{enumerate}
    \item $\Ks \Ks^{-1}=$ Id,
    \item $\Ks^{-1}(b,w)\rightarrow 0$, when $|b-w|\rightarrow\infty$.
\end{enumerate}
The rational functions $\fs_{wx}(z)$ are the analog of the rational
functions $f_{wv}(z)$ of \cite{Kenyon3}, and are defined in the
following way. Let $w\in W$, and let
$x\in B$ (resp. $x\in W$); consider the edge-path
$w=w_1,b_1,\ldots,w_k, b_k=x$ (resp.
$w=w_1,b_1,\ldots,w_k,b_k,w_{k+1}=x$) of $G^*$ from $w$ to $x$.
Let $R(w_j b_j)$ be the rhombus associated to the edge $w_j b_j$,
and denote by $w_j,x_j,b_j,y_j$ its vertices in cclw
(counterclockwise) order;
$e^{i\alpha_j}$ is the complex vector $y_j-w_j$, and
$e^{i\beta_j}$ is the complex vector $x_j-w_j$. In a similar way,
denote by $w_{j+1},x_j',b_j,y_j'$ the vertices of the rhombus
$R(w_{j+1}b_j)$ in cclw order, then $e^{i\alpha_j'}$ is the
complex vector $y_j'-w_{j+1}$, and $e^{i\beta_j'}$ is the complex
vector $x_j'-w_{j+1}$. The function $\fs_{wx}(z)$ is defined inductively along
the path,
\begin{eqnarray*}\label{61}
\fs_{ww}(z)&=&1,\\
\fs_{w b_j}(z)&=&\fs_{w w_j}(z)
\frac{(-1)^{\II_{(w_j,b_j)}}e^{i\frac{\alpha_j+\beta_j}{2}}}{(z-e^{i\alpha_j})(z-e^{i\beta_j})},\\
\fs_{w w_{j+1}}(z)&=&\fs_{w b_j}
(z)(-1)^{\II_{(w_{j+1},b_j)}}e^{-i\frac{\alpha_j'+\beta_j'}{2}}(z-e^{i\alpha_j'})(z-e^{i\beta_j'}).\\
\end{eqnarray*}
By Lemma $9$ of \cite{Bea}, the function $\fs_{wx}$ is well
defined. Moreover, Lemma $10$ of \cite{Bea} gives an explicit
expression for $\Ks^{-1}$ (it is the analog of Theorem $4.2$ of \cite{Kenyon3}):
\begin{lem}{\rm\bf\cite{Bea}}\label{prop61}
The inverse real Dirac operator is given by
\begin{equation}\label{68}
\Ks^{-1}(b,w)=\frac{1}{4\pi^2 i}\int_{C} \fs_{wb}(z) \log z \;dz,
\end{equation}
where $C$ is a closed contour surrounding cclw the part of the
circle $\{ e^{i \theta} \,|\, \theta \in [\theta_0 - \pi +
\Delta,\theta_0 + \pi - \Delta ]\}$, which contains all the poles
of $\fs_{wb}$, and with the origin in its exterior.
\end{lem}

\subsection{Partition function of the torus}\label{subsec12}

\noindent In this section, we give the explicit formula of Kasteleyn
\cite{Kasteleyn2} (in the $\ZZ^2$ case), and Tesler \cite{Tesler} (in
the general periodic case) for the toroidal partition function
$Z(G_n,\nu)$. Note that this result was already described in
\cite{Bea}, we nevertheless choose to repeat it here since its
understanding is important for the proof of Theorem \ref{thm1}.\jump
Let us first orient the edges of $G$. Consider the
graph $G_1$, it is a bipartite graph on the torus. Fix a
reference matching $M_0$ of $G_1$. For every other perfect
matching $M$ of $G_1$, consider the superposition $M\cup M_0$
of $M$ and $M_0$, then $M\cup M_0$ consists of doubled edges and
cycles. Let us define four parity classes for perfect matchings
$M$ of $G_1$: (e,e) consists of perfect matchings $M$, for
which cycles of $M\cup M_0$ circle the torus an even number of
times horizontally and vertically; (e,o) consists of perfect
matchings $M$, for which cycles of $M\cup M_0$ circle the torus
and even number of times horizontally, and an odd number of times
vertically; (o,e) and (o,o) are defined in a similar way. By
Tesler \cite{Tesler}, one can construct an orientation of the
edges of $G_1$, so that the corresponding adjacency matrix
$\Ks_1^1$ has the following property: perfect matchings which
belong to the same parity class have the same sign in the
expansion of the determinant; moreover of the four parity classes,
three have the same sign and one the opposite sign. By an
appropriate choice of sign, we can make the (e,e) class have the
plus sign in $\det \Ks_1^1$, and the other three have minus sign.
Consider a horizontal and a vertical cycle of the dual graph $G_1^*$
of $G_1$. Then
define $\Ks_2^1$ (resp. $\Ks_3^1$) to be the matrix $\Ks_1^1$
where the sign of the coefficients corresponding to edges crossing
the horizontal (resp. vertical) cycle is reversed; and define
$\Ks_4^1$ to be the matrix $\Ks_1^1$ where the sign of the
coefficients corresponding to the edges crossing both cycles are
reversed. By Kasteleyn \cite{Kasteleyn2} (in the $\ZZ^2$ case), and
Tesler \cite{Tesler} (in the general case), we have the following,
\begin{equation*}\label{616}
Z(G_1,\nu)=\frac{1}{2}(-\det \Ks_1^1+\det \Ks_2^1+\det
\Ks_3^1+\det \Ks_4^1).
\end{equation*}
The orientation of the edges of $G_1$ defines a periodic
orientation of the graph $G$. For every $n$, consider the graph
$G_n$, and the four matrices
$\Ks_1^n,\Ks_2^n,\Ks_3^n,\Ks_4^n$ defined as above. These matrices
are called the {\bf Kasteleyn matrices} of the graph $G_n$.
\begin{thm}{\rm \cite{Kasteleyn2,Tesler}}\label{thm61}
\begin{equation*}\label{617}
Z(G_n,\nu)=\frac{1}{2}(-\det\Ks_1^n+\det\Ks_2^n+\det\Ks_3^n+\det\Ks_4^n).
\end{equation*}
\end{thm}
The orientation defined on the edges of the graph $G$ is a
clockwise odd orientation. Let $\Ks$ be the real Dirac operator
indexed by the vertices of $G$, corresponding to 
this clockwise odd orientation. Note that except for edges crossing
the horizontal and the vertical
cycle, the coefficients of the Kasteleyn matrices $\Ks_\l^n$ and
of the real Dirac operator agree on edges they have in common.

\subsection{Proof}\label{subsec13}

\noindent Consider an orientation of the edges of $G$ defined
as in Section \ref{subsec12}, and let $\Ks$ be the real Dirac
operator indexed by vertices of $G$, corresponding to this
orientation. Let $\Ks_1^n,\Ks_2^n,\Ks_3^n,\Ks_4^n$ be the
Kasteleyn matrices of the graph $G_n$. The following computations are
inspired form \cite{Kenyon3}. By linear algebra, for
every $\l=1,\ldots,4$, and for every edge $w_i b_i$ of $G_n$,
we have
\begin{equation*}\label{335}
\frac{\partial\det\Ks_\l^n}{\partial
(\Ks_\l^n(w_i,b_i))}=(\det\Ks_\l^n)((\Ks_\l^n)^{-1}(b_i,w_i)).
\end{equation*}
Denote by $e_1=w_1 b_1,\ldots,e_m=w_m b_m$ the edges of
$G_1$, and let $\theta_1,\ldots,\theta_m$ be the
corresponding rhombus angles. Since the graph $G_n$ is
invariant under $\Lambda$-translates, we know that for every
$\Lambda$-translate $w_i^t b_i^t$ of the edge $w_i b_i$, the
coefficient $\Ks_\l^n(w_i^t,b_i^t)=\pm \Ks_\l^n(w_i,b_i)$. The
minus sign only occurs when $\l=2,3,4,$ (recall that in the
definition of $\Ks_2^n,\ldots,\Ks_4^n$, the sign of the entries
which cross the horizontal and/or the vertical cycle of $G_n^*$
is reversed), but then
$(\Ks_\l^n)^{-1}(b_i^t,w_i^t)=\pm(\Ks_\l^n)^{-1}(b_i,w_i)$. Hence,
for every $\l=1,\ldots,4$, and for every $i=1,\ldots,m$,
\begin{eqnarray}\label{318}
\frac{\partial \det \Ks_\l^n}{\partial\theta_i}&=&
\sum_{\mbox{{\scriptsize $w_i^t b_i^t$ translates of $w_i b_i$}}}
\frac{\partial\det\Ks_\l^n}{\partial(\Ks_\l^n(w_i^t,b_i^t))}
\frac{\partial(\Ks_\l^n(w_i^t,b_i^t))}{\partial\theta_i},\nonumber\\
&=&n^2(\det\Ks_\l^n)((\Ks_\l^n)^{-1}(b_i,w_i))\frac{\partial
\Ks_\l^n(w_i,b_i)}{\partial\theta_i},
\end{eqnarray}
where the edges $e_1,\ldots,e_m$ do not cross the horizontal and
the vertical path of the dual graph $G_n^*$ of $G_n$. Define the function
$\vphi_n(\theta_1,\ldots,\theta_m)$ by
\begin{equation*}
\vphi_n(\theta_1,\ldots,\theta_m)=\frac{1}{n^2}\log
Z(G_n,\nu).
\end{equation*}
By Theorem \ref{thm61}, we have
$Z(G_n,\nu)=\frac{1}{2}(-\det\Ks_1^n+\det\Ks_2^n+\det\Ks_3^n+\det\Ks_4^n)$.
Moreover, for every $\l=1,\ldots,4$,
$\Ks_\l^n(w_i,b_i)=\Ks(w_i,b_i)$, so that using (\ref{318}), we
obtain for every $i=1,\ldots,m$, {\scriptsize
\begin{equation}\label{324}
\frac{\partial
\vphi_n}{\partial\theta_i}(\theta_1,\ldots,\theta_m)=
\left(\frac{\partial\Ks(w_i,b_i)}{\partial\theta_i}\right)
\left(\frac{-\det\Ks_1^n}{2Z(G_n,\nu)}(\Ks_1^n)^{-1}(b_i,w_i)+\sum_{\l=2}^4
\frac{\det\Ks_\l^n}{2Z(G_n,\nu)}(\Ks_\l^n)^{-1}(b_i,w_i)
\right).
\end{equation}}
The next part of the argument can be found in \cite{Kenyon0}. The
second bracket of equation (\ref{324}) is a weighted average of
the four quantities $(\Ks_\l^n)^{-1}(b_i,w_i)$, with weights
$\pm\frac{1}{2}\det\Ks_\l^n/Z(G_n,\nu)$. These weights are
all in the interval $(-1/2,1/2)$ since for every $\l=1,\ldots,4$,
$2Z(G_n,\nu)>|\det\Ks_\l^n|$. Indeed, $Z(G_n,\nu)$
counts the weighted sum of dimer configurations of $G_n$,
whereas $|\det\Ks_\l^n|$ counts some configurations with negative
sign. Since the weights sum to $1$, the weighted average converges
to the same value as each $(\Ks_\l^n)^{-1}(b_i,w_i)$. 

By Proposition $2$ of \cite{Bea}, for every $\l=1,\ldots,4$,
$(\Ks_\l^n)^{-1}(b_i,w_i)$ converges to $\Ks^{-1}(b_i,w_i)$ on a
subsequence $(n_j)$ of $n$'s. Hence, for every $\eps>0$, there
exists $n_0$ such that for $n\geq n_0$, $n\in(n_j)$, equation
(\ref{324}) can be written as
\begin{equation}\label{equ4}
\frac{\partial
\vphi_n}{\partial\theta_i}(\theta_1,\ldots,\theta_m)=
\frac{\partial\Ks(w_i,b_i)}{\partial\theta_i}\Ks^{-1}(b_i,w_i)+\eps.
\end{equation}
Let us compute the right hand side of (\ref{equ4}), using the
notations of Section \ref{subsubsec112}. By definition,
$\Ks(w_i,b_i)=(-1)^{\II_{(w_i,b_i)}}2\sin\theta_i$, hence
\begin{equation}\label{equ5}
\frac{\partial\Ks(w_i,b_i)}{\partial\theta_i}=(-1)^{\II_{(w_i,b_i)}}2\cos\theta_i.
\end{equation}
Moreover, by definition $\fs_{w_i
  b_i}(z)=(-1)^{\II_{(w_i,b_i)}}\frac{e^{i\frac{(\alpha_i+\beta_i)}{2}}}{(z-e^{i\alpha_i})(z-e^{i\beta_i})}$. So that using Lemma \ref{prop61}, and the Residue Theorem yields:
\begin{eqnarray}\label{equ6}
\Ks^{-1}(b_i,w_i)&=&
  \frac{1}{4\pi^2 i}(-1)^{\II_{(w_i,b_i)}}e^{i\frac{(\alpha_i+\beta_i)}{2}}
  \int_{C}\frac{\log z}{(z-e^{i\alpha_i})(z-e^{i\beta_i})}dz,\nonumber\\
&=&\frac{1}{2\pi}(-1)^{\II_{(w_i,b_i)}}e^{i\frac{(\alpha_i+\beta_i)}{2}}
   \left(\frac{i\alpha_i}{e^{i\alpha_i}-e^{i\beta_i}}+\frac{i\beta_i}{e^{i\beta_i-\alpha_i}}\right),\nonumber\\
&=&\frac{1}{2\pi}(-1)^{\II_{(w_i,b_i)}}\frac{\theta_i}{\sin\theta_i}.
\end{eqnarray}
Combining equations (\ref{equ5}) and (\ref{equ6}), implies:
\begin{equation*}
\frac{\partial
\vphi_n}{\partial\theta_i}(\theta_1,\ldots,\theta_m)=
\frac{\theta_i}{\pi}\mbox{cotan}\,\theta_i+\eps.
\end{equation*}
By \cite{Kenyon3}, there is a continuous way to deform the graph
$G$ so that all rhombus angles tend to $0$ or $\pi/2$. The same
transformation can be applied to $G_n$, for every $n$. Denote
by $\theta_i^0$ the angle $\theta_i$ after such a deformation. Let
$M$ be the number of angles $\theta_1^0,\ldots,\theta_m^0$ which
are equal to $\pi/2$. By the argument of \cite{Kenyon3}, for
$n\geq n_0$, $n\in(n_j)$, we obtain {\scriptsize
\begin{equation}\label{345}
\vphi_n(\theta_1,\ldots,\theta_m)=\sum_{i=1}^m\left(
\frac{\theta_i}{\pi}\log
2\sin\theta_i+\frac{1}{\pi}L(\theta_i)\right)
-\frac{M}{2}\log2
+\vphi_n(\theta_1^0,\ldots,\theta_m^0)+\eps.
\end{equation}}

\noindent Let us compute $\vphi_n(\theta_1^0,\ldots,\theta_m^0)$.
Suppose the above deformation is applied to the graph $G$, then
edges corresponding to rhombus angles $0$ have weight $0$, and
those corresponding to rhombus angles $\pi/2$ have weight $2$.
Removing the $0$ weight edges, the deformed graph consists of
independent copies of $\ZZ$. Applying the deformation to
$G_1$, we obtain graphs $\ZZ/m_1\ZZ,\ldots,\ZZ/m_p\ZZ$, for
even $m_1,\ldots,m_p$, where each $\ZZ/m_j\ZZ$ consists of weight
$2$ edges. We can compute,
\begin{equation*}\label{346}
Z(G_n,\nu)=(2^{\frac{n m_1}{2}}2)^n\ldots (2^{\frac{n
m_p}{2}}2)^n=2^{n^2\frac{M}{2}}2^{np}.
\end{equation*}
Hence, for every $\eps>0$, there exists $n_1$, such that for every
$n\geq n_1$, we have
\begin{equation}\label{347}
\vphi_n(\theta_1^0,\ldots,\theta_m^0)=\frac{M}{2}\log
2+\eps.
\end{equation}
Combining equations (\ref{345}) and (\ref{347}) yields, for every
$n\geq \max\{n_0,n_1\}$, $n\in (n_j)$,
\begin{equation*}\label{348}
\vphi_n(\theta_1,\ldots,\theta_m)=
\sum_{i=1}^m\left(\frac{\theta_i}{\pi}\log
2\sin\theta_i+\frac{1}{\pi}L(\theta_i)\right)+ 2\eps.
\end{equation*}
This implies,
\begin{equation*}\label{349}
{\lim}'_{\,n\rightarrow\infty}\vphi_n(\theta_1,\ldots,\theta_m)=\sum_{i=1}^m
\left(\frac{\theta_i}{\pi}\log
2\sin\theta_i+\frac{1}{\pi}L(\theta_i)\right).
\end{equation*}
where the limit is taken on the subsequence $(n_j)$. Using Theorem
3.5 of \cite{KOS}, we deduce that
$\vphi_n(\theta_1,\ldots,\theta_m)$ converges (as
$n\rightarrow\infty$) to the same limit as the above subsequence.
\hspace*{\fill}$\square$

\bibliographystyle{alpha}

\end{document}